\documentclass[11pt]{article}
\usepackage{amssymb,amsfonts,amsmath,amsthm,cite}
\usepackage{graphicx,epsfig,subfigure}
\usepackage{enumerate}
\usepackage{tikz}

\newtheorem{thm}{Theorem}[section]

\newtheorem{lem}{Lemma}[section]

\newtheorem{prob}{Problem}[section]

\makeatletter \@addtoreset{equation}{section}
\def\red{}
\def\pf{\noindent {\it Proof.\ }}
\def\qed{\hfill \rule{4pt}{7pt}}

\title{\bf  The Tur\'{a}n Number for Spanning Linear Forests\thanks
{The work was supported by NNSF of China (No.11701407, No.11671296, No.61502330).}}
\author{Jian Wang\footnote{\ Corresponding author:
wangjian01@tyut.edu.cn}, Weihua Yang}
\date{\small
{\it Department of Mathematics,\\
\vskip4pt Taiyuan University of Technology, Taiyuan, 030024, P.R.China}
}

\begin{document}

\maketitle

\noindent {\bf Abstract.}   {For a set of graphs $\mathcal{F}$}, the extremal number {$ex(n;\mathcal{F})$}
is the maximum number of edges in a graph of order $n$ not containing {any} subgraph isomorphic to
some graph in {$\mathcal{F}$}. If {$\mathcal{F}$} contains a graph on $n$ vertices, then we often call the problem {a} spanning Tur\'{a}n problem. A linear forest is a graph whose connected components are all paths and isolated vertices. In this paper, we {let} $\mathcal{L}_n^k$ be the set of all linear forests of order $n$ with at least $n-k+1$ edges. {We} prove that when $n\geq 3k$ \red{and $k\geq 2$},
\[
ex(n;\mathcal{L}_n^k)=\binom{n-k+1}{2}+ O(k^2).
\]
Clearly, the result is {interesting} when $k=o(n)$.


\allowdisplaybreaks

\section{Introduction}
\label{sec-intro}
\red{For {a set of graphs $\mathcal{F}$}, the extremal number {$ex(n;\mathcal{F})$}
is the maximum number of edges in a graph of order $n$ not containing {any} subgraph isomorphic to some graph in {$\mathcal{F}$}. {Tur\'{a}n introduced this problem in \cite{turan}, and we recommend \cite{sidorenko,keevash} for surveys on Tur\'{a}n problems for graphs and hypergraphs.} {Ore}\cite{ore} proved that a non-Hamiltonian graph of order $n$ has at most $\binom{n-1}{2}+1$ edges. Ore's theorem can also be expressed as $ex(n;\{C_n\})=\binom{n-1}{2}+1$ where $C_n$ is the cycle of order $n$. Similarly, $ex(n;\{P_n\})=\binom{n-1}{2}$ where $P_n$ is the {path} of order $n$.}

Recently, many {generalizations} of Ore's theorem have been studied. {Alon and Yuster}\cite{alon} \red{extended} Ore's result to spanning structures other than just Hamilton cycles. They prove that if $H$ is a graph of order $n$ with minimum degree $\delta(H) > 0$, and  maximum degree $\Delta(H)\leq \sqrt{n}/40$, then $ex(n;\{H\})=\binom{n-1}{2}+\delta(H)-1$ assuming $n$ is sufficiently large. Moreover, some researchers extend Ore's theorem to the hypergraph {setting}\cite{katona,tuza,glebov}.

Define the {\it Hamiltonian completion number} of a graph $G$, denoted by $h(G)$, to be the minimum
number of edges that need be added to make $G$ Hamiltonian. The Hamiltonian completion problem was introduced in \red{the} 1970s by Goodman and Hedetniemi\cite{good74,good75}. From {the} viewpoint of Hamiltonian completion number of graphs, Ore's Theorem can also be restated as a graph with {$h(G)>0$} of order $n$ has at most $\binom{n-1}{2}+1$ edges. Then, it's natural to try to extend Ore's result to {graphs} with $h(G)\geq k$ for some $k\geq 1$.  Let $\mathcal{L}_n^k$ be the set of all linear forests of order $n$ with at least $n-k+1$ edges. Clearly, the problem is equivalent to {determining} the {Tur\'{a}n} number $ex(n;\mathcal{L}_n^k)$. In this paper, we prove that {when $n\geq 3k$ and $k\geq 2$},
\[
ex(n;\mathcal{L}_n^k)=\binom{n-k+1}{2}+O(k^2).
\]

{It should be mentioned that Lidick\'{y} et. al. determine the Tur\'{a}n number of linear forests of arbitrary order for $n$ sufficiently large in \cite{lidicky}. However, the number of vertices in their forbidden linear forests does not depend on $n$.}
The rest of this short paper is organized as follows. In Section 2, we give lower bounds on {the} Tur\'{a}n number $ex(n;\mathcal{L}_n^k)$ and prove a useful lemma. In Section 3, we prove the main theorem. In Section 4, we give some concluding remarks.

\section{Lower Bounds and A Useful Lemma}

Let $G_0$ be the union graph of $K_{n-k+1}$ and $k-1$ isolated vertices. It is easy to see that $G_0$ does not contain any linear forest with more than $n-k$ edges. Thus, we have the following lemma.

\begin{lem}\label{lowerbounds}
The Tur\'{a}n number $ex(n;\mathcal{L}_n^k)$ has the following lower bound.
\[
ex(n;\mathcal{L}_n^k)\geq e(G_0) = \binom{n-k+1}{2}.
\]
\end{lem}

\begin{lem}[\red{\hspace{1sp}\cite{bondy}}]\label{ore}
\red{Let $G$ be a graph with $n$ vertices and let $P$ be} a Hamiltonian path in $G$ with {endpoints} $u$ and $v$. If {$d(u)+d(v)\geq n$}, then $G$ contains a Hamiltonian cycle.
\end{lem}

For simplicity, we view isolated vertices as paths of length zero, whose end vertices are the same.
\begin{lem}\label{lfedges}
{Let $k\geq 2$ be an integer, and suppose that $G$ is a graph that contains a spanning linear forest $F$ with $n-k$ edges. If $u$ and $v$ are vertices that are endpoints of different paths in $F$ and $d(u)+d(v)\geq n-k+1$, then $G$ contains a spanning linear forest with $n-k+1$ edges.}
\end{lem}

\pf
Let $H$ be the join graph of $G$ and $k-1$ isolated vertices $v_1,v_2,\ldots,v_{k-1}$. That means
\[V(H)=V(G)\cup \{v_1,v_2,\ldots,v_{k-1}\}\]
and
\[ E(H)=E(G)\cup \{v_iw\colon i=1,2,\ldots, k-1 \mbox{ and } w\in V(G)\}.\]
 {Let $P_1,P_2,\ldots,P_k$ be $k$ paths of $F$, and assume that $u$ is the first vertex of $P_1$, and $v$ is the last vertex of $P_k$.} Then $P_1v_1P_2v_2\ldots,v_{k-1}P_k$ forms a Hamiltonian path
of $H$ with endpoints $u$ and $v$. Moreover, $d_{H}(u)+d_{H}(v)\geq n-k+1+2(k-1)=n+k-1$. By Lemma \ref{ore}, it follows that
$H$ contains a Hamiltonian cycle. Finally, by removing vertices $v_1,v_2,\ldots,v_{k-1}$ from this Hamiltonian cycle, we obtain a linear forest of $G$ with $n+k-1-2(k-1)=n-k+1$ edges. Thus, the lemma holds. \qed

\section{The Main Result}

\begin{thm}
For $n\geq 3k$ \red{and $k\geq 2$}, the Tur\'{a}n number $ex(n;\mathcal{L}_n^k)$ has the following upper bound.
\[
ex(n;\mathcal{L}_n^k)\leq \binom{n-k+1}{2}+\frac{k^2-3k+4}{2}.
\]
\end{thm}
\pf
{Let $G(V,E)$ be any graph with $n$ vertices, and the maximum number of edges subject to $h(G)\geq k$. We may assume that $h(G)=k$ for if  $h(G)>k$, then one may add an edge to $G$ to obtain a new graph $G'$ with $h(G')\geq k$, $e(G')>e(G)$ and having the same number of vertices as $G$. As $h(G)=k$, we can choose a linear forest $F$ with $n-k$ edges and furthermore, among all such subgraphs, we may choose $F$ so that it has the fewest number of isolated vertices.}

Since $F$ has $n$ vertices and $n-k$ {edges, then} $F$ has $k$ connected components {consisting} of paths and isolated vertices.
%Suppose we have $i$ isolated vertices and $k-i$ paths.
By Lemma \ref{lfedges}, for any two end vertices $u, v$ in different {connected} components {of $F$}, we have $d(u)+d(v)\leq n-k$. Otherwise, $G$ contains a linear forest with $n-k+1$ edges, {contradicting with $h(G)\geq k$}.

We claim that all but at most two end vertices of $F$ have degree at most $\frac{n-k}{2}$. Moreover, if two end vertices of $F$ have degree greater than $\frac{n-k}{2}$, {then} they \red{are end vertices of the same path}. If there are three vertices with degree greater than $\frac{n-k}{2}$, then at {least} two of them fall into the different components, assume they are $u$ and $v$. Then $d(u)+d(v)> n-k$. {It follows that $G$ contains a linear forest with $n-k+1$ edges, contradicting with $h(G)\geq k$}. Thus, the claim holds.

{The remainder of the proof splits into two cases, depending on whether or not $F$ contains isolated vertices. For each case, we give the upper bound on the number of edges in $G$.}

{{\bf Case 1.}\,  There {are} no isolated vertices in $F$.}   {Then we claim that} any two end vertices in different {paths} of $F$ are not adjacent. Otherwise, by {using the} edge between these two vertices we obtain a linear forest with $n-k+1$ edges, {contradicting with $h(G)\geq k$}. Thus, each end vertex of $F$ has degree at most $n-1-(2k-2)=n-2k+1$.
Let $X$ be the set of all the end vertices of \red{the} $k$ paths in $F$ and $G'=G[V\backslash X]$. Since $G'$ has $n-2k$ vertices, we have
\[
e(G')\leq \binom{n-2k}{2}.
\]
{Now we bound the number of edges of $G$ as follows:}
 \begin{align}\label{noislated}
e(G) &\leq  \binom{n-2k}{2}+2(n-2k+1)+2(k-1)\cdot\frac{n-k}{2} \nonumber\\[5pt]
&=\binom{n-k+1}{2}+\frac{k^2-3k+4}{2}.
\end{align}

{{\bf Case 2.}\, There {are} $i$ isolated vertices in $F$ for some $i$ with $0<i<k$.} {Let $X=\{x_1,x_2,\ldots,x_i\}$ be} the set of $i$ isolated vertices in $F$. We claim that any vertex in $X$ has degree at most $k-i-1$. Let $x_s\in X$ and $P=y_1y_2\ldots y_m$ be any path in $F$. If $x_sy_t$ is an edge in $G$ for some $1<t<m$, then $m$ has to be $3$. Otherwise, if $t=2$, then by replacing $P$ and $x_s$ with $y_1y_2x_s$ and $y_3\ldots y_m$ from $F$, we obtain a linear forest with less isolated vertices, which contradicts with the selection of $F$. If $t\geq 3$, {then} by replacing $P$ and $x_s$ with $y_1\ldots y_{t-1}$ and $x_sy_t\ldots y_m$ from $F$, we obtain a linear forest with less isolated vertices, a contradiction. Therefore, neighbors of $x_s$ can only be internal vertices of paths of length \red{two} in $F$. Since $n\geq 3k$, $F$ contains at most $k-i-1$ paths of length \red{two}. Thus, any vertex in $X$ has degree at most $k-i-1$ in {$G$}.

Let $G'=G[V\backslash X]$ and $F'=F[V\backslash X]$. We claim that $F'$ {is a} spanning subgraph of $G'$ that is a linear forest with maximum number of edges. Otherwise, if $G'$ has a linear forest $F^*$ with more edges\red{, then} by replacing edges in $F'$ with those in $F^*$, we obtain a spanning linear forest of $G$ with more edges, a contradiction. Since $F'$ has $n-i$ vertices and $k-i$ components and contains no isolated vertices, by inequality {\eqref{noislated}} we have
 \begin{align*}
e(G')& \leq \binom{(n-i)-(k-i)+1}{2}+\frac{(k-i)^2-3(k-i)+4}{2} \\[5pt]
&= \binom{n-k+1}{2}+\frac{(k-i)^2-3(k-i)+4}{2}.
\end{align*}
 Note that $1\leq i\leq k-1$, {we bound the number of edges of $G$ as follows:}
 \begin{align*}
e(G)&\leq  e(G')+i(k-i-1) \\[5pt]
&\leq \binom{n-k+1}{2}+\frac{(k-i)^2-3(k-i)+4}{2}+i(k-i-1)\\[5pt]
&= \binom{n-k+1}{2}+\frac{k^2-3k}{2}-\frac{1}{2}\cdot\left(i-\frac{1}{2}\right)^2+\frac{17}{8}\\[5pt]
&\red{\leq} \binom{n-k+1}{2}+\frac{k^2-3k+4}{2}.
\end{align*}
{Combining} the two cases, we conclude that
\[ex(n;\mathcal{L}_n^k)\leq \binom{n-k+1}{2}+\frac{k^2-3k+4}{2}.\]
\qed

\section{Concluding Remarks}

\red{In this paper, we prove that {when $n\geq 3k$ and $k\geq 2$},
\[
\binom{n-k+1}{2} \leq ex(n;\mathcal{L}_n^k)\leq \binom{n-k+1}{2}+\frac{k^2-3k+4}{2},
\]
the result is {interesting} when $k=o(n)$. For $k=2$,  $\mathcal{L}_n^k$ denote the set of all linear forests of order $n$ with at least $n-1$ edges, which is exactly the set of Hamiltonian paths. Thus, the lower bound is reachable for $k=2$.  Furthermore, we guess that there exists a constant $k_0$ such that $ex(n;\mathcal{L}_n^k)=\binom{n-k+1}{2}$ for $k<k_0$. And we end up this paper by proposing the following two problems.}

\begin{prob}
Determine the exact value of $ex(n;\mathcal{L}_n^k)$ for $k=o(n)$.
\end{prob}
\begin{prob}
Let $c$ be a constant satisfying {$0<c<1$. Determine} the value of $ex(n;\mathcal{L}_n^k)$ for $k=cn$.
\end{prob}

 \noindent {\bf Acknowledgments.}
The authors would like to express their gratitude to the anonymous referees for their kind and detailed comments on the original manuscript, which improve the presentation of the manuscript substantially.

\end{document}